\documentclass[reqno]{amsart}
\newif\ifscreenmode
\newif\ifarXivmode
\screenmodetrue   
\arXivmodetrue    

\ifscreenmode
\newlength{\screenmodehorizontalmargin}
\newlength{\screenmodeverticalmargin}
\ifarXivmode
\fi
\usepackage[paperwidth={\dimexpr 360pt + 2\screenmodehorizontalmargin\relax},
            paperheight={\dimexpr 584pt + 34pt + 2\screenmodeverticalmargin\relax},
            textwidth=360pt,
            textheight=584pt,
            hmarginratio=1:1,
            vmarginratio=1:1,
            includeheadfoot]{geometry}
\fi

\usepackage[english]{babel}
\usepackage{amsmath,amssymb,mathtools}
\usepackage{amsthm,thmtools}
\usepackage[shortlabels]{enumitem}
\usepackage{lmodern,microtype}
\usepackage{tikz}

\usepackage[colorlinks=true,linkcolor=,urlcolor=,citecolor=]{hyperref} 
\usepackage[capitalize,nameinlink]{cleveref} 

\hypersetup{hypertexnames=false} 

\usetikzlibrary{cd}
\usetikzlibrary{positioning,calc,arrows.meta}
\usetikzlibrary{decorations.pathmorphing,decorations.pathreplacing,patterns,shapes.symbols}

\newcommand{\ang}[1]{\langle #1 \rangle} 
\newcommand{\wtd}{\widetilde}

\DeclareRobustCommand{\SkipTocEntry}[5]{}

\allowdisplaybreaks

\setlist{itemsep=.5\baselineskip,topsep=.5\baselineskip}

\numberwithin{equation}{section}
\declaretheorem[style=plain,numberlike=equation]{theorem}
\declaretheorem[style=plain,numberlike=theorem]{proposition}
\declaretheorem[style=plain,numberlike=theorem]{lemma}
\declaretheorem[style=plain,numberlike=theorem]{corollary}

\declaretheorem[style=definition,numberlike=theorem]{definition}
\declaretheorem[style=definition,numberlike=theorem]{construction}
\declaretheorem[style=definition,numberlike=theorem,title=Definition]{defn}

\declaretheorem[style=remark,numberlike=theorem]{example}

\newcommand{\iso}{\cong}
\newcommand{\incl}{\hookrightarrow}
\DeclareMathOperator{\vspan}{span}

\newcommand{\C}{\mathbb{C}}
\newcommand{\Z}{\mathbb{Z}}

\newcommand{\Q}{\mathbb{Q}}
\newcommand{\K}{\mathbb{K}}

\newcommand{\mcF}{\mathcal{F}}
\newcommand{\mcG}{\mathcal{G}}

\newcommand{\mcN}{\mathcal{N}}

\title[Adian--Rabin theorem for quantum automorphism groups]{An Adian--Rabin theorem for quantum automorphism groups of graphs}

\author[J. \MakeLowercase{\textsc{van}} Dobben de Bruyn]{Josse van Dobben de Bruyn}
\address{Charles University, Department of Applied Mathematics, Malostranské náměstí 25, 118~00 Praha~1, Czech Republic}
\email{josse.van-dobben-de-bruyn@matfyz.cuni.cz}

\author[W. Slofstra]{William Slofstra}
\address{University of Waterloo, Department of Pure Mathematics and Institute for Quantum Computing, 200 University Avenue West, Waterloo, ON, Canada N2L 3G1}
\email{weslofst@uwaterloo.ca}


\newcommand{\MP}{{\operatorname{MP}}}
\newcommand{\simple}{{\operatorname{sim}}}
\DeclareMathOperator{\Qut}{Qut}
\DeclareMathOperator{\SO}{SO}
\DeclareMathOperator{\Aut}{Aut}
\newcommand{\FP}{\mathtt{FP}}
\newcommand{\CCQG}{\mathtt{CCQG}}
\newcommand{\Ast}{*}

\newcommand{\GrphInf}{\mcG_{\operatorname{inf}}}
\newcommand{\GrphTriv}{\mcG_{\operatorname{triv}}}
\newcommand{\AR}{\mathcal{K}}

\DeclareMathOperator{\id}{id}
\newcommand{\tensor}{\otimes}
\DeclareMathOperator{\Iso}{Iso}

\begin{document}

\begin{abstract}
    Given a graph $G$, does it have quantum symmetries, or in other words, is
    the algebra of the quantum automorphism group non-commutative? In the world of ordinary
    groups, determining whether a finitely presented group is commutative is
    undecidable. This is an immediate consequence of the Adian--Rabin theorem,
    which shows that many natural properties of finitely presented groups are
    undecidable. In this paper, we prove an Adian--Rabin theorem for quantum
    automorphism groups of graphs, showing that many natural properties of the
    quantum automorphism group of a graph (including having quantum symmetries)
    are undecidable. We also show that it's undecidable to determine if the 
    discrete quantum automorphism group is non-trivial. Our methods are based on
    solution groups of linear system nonlocal games, and we also prove an
    Adian--Rabin theorem in this context. 
\end{abstract}

\maketitle

\section{Introduction}

The quantum automorphism group $\Qut(X)$ of a finite graph $X$ is a compact
quantum group which captures the noncommutative symmetries of the graph.
We use the Hopf-algebraic framework, where a compact quantum group is a Hopf $*$-algebra with a Haar state.
The abelianization of $\Qut(X)$ is the algebra $C(\Aut(X))$ of functions
on the classical automorphism group $\Aut(X)$. If $\Qut(X)$ is noncommutative,
then $X$ is said to have quantum symmetries. There are many examples of graphs
with quantum symmetries, and indeed, it's well-known that $\Qut(X)$ can be
infinite dimensional (e.g. \cite[Example~1.1.8]{Neshveyev-Tuset}), whereas
$\Aut(X)$ is finite by definition.
A recent result of the first author and others is that there exist graphs $X$ where
$\Qut(X)$ is non-trivial and $\Aut(X)$ is trivial \cite{Dobben-Roberson-Schmidt}.
For most graphs, $\Qut(X)$ is trivial \cite[Theorem~3.14]{Lupini-Mancinska-Roberson}.

It's natural to ask whether there is a procedure to determine, for a given graph
$X$, whether $\Qut(X)$ is trivial or has one of the other properties mentioned
above. Famously, the same type of question comes up with finitely presented
groups: given a finitely presented group, we'd like to know if the group is
trivial, abelian, finite, etc. For groups, the Adian--Rabin theorem \cite{Adyan,Rabin} states that
if (P) is a Markov property, then determining whether a finitely presented
group has (P) is undecidable. A Markov property is a property of groups (P)
such that there is a finitely presented group $G_+$ which has (P), and another
finitely presented group $G_-$ such that any finitely presented group containing $G_-$ as a
subgroup does not have (P). All of the group properties just listed are Markov
properties, and hence are undecidable.

The first main result of this paper is an Adian--Rabin theorem for quantum
automorphism groups:
\begin{theorem}
	\label{thm:main}
	Let $(P)$ be a property of quantum automorphism groups of graphs, for which there exists a graph
	$Y_+$ such that $\Qut(Y_+)$ has property $(P)$, and a graph $Y_-$ such that
	$\Qut(Y_-)$ is cocommutative and is not a Hopf $*$-subalgebra 
    of any $\Qut(X)$ with property $(P)$.  Then it is undecidable (RE-hard) to
    determine for an input graph $X$ whether or not $\Qut(X)$ has property $(P)$.
\end{theorem}
Being commutative (as a Hopf $*$-algebra), being trivial, and being
finite-dimensional are all properties which satisfy the hypotheses of
\Cref{thm:main}, and hence are undecidable for quantum automorphism groups. We
do not know whether we can drop the requirement that $\Qut(Y_-)$ is
cocommutative. If $\Qut(Y_+)$ is cocommutative as well, then the condition on
$\Qut(Y_-)$ can be slightly weakened: we only need to check that it is not a
Hopf $*$-subalgebra of any cocommutative $\Qut(X)$ with property $(P)$. This
variation, along with the proof of \Cref{thm:main}, is in \Cref{sec:qag}.

The proof of the Adian--Rabin theorem for groups uses an elementary construction:
Given a finitely presented group $G = \ang{S : R}$ and a word $w \in \mcF_S$ in the
free group generated by $S$, it's possible to add generators and relations to
get a new finitely presented group $\wtd{G}_w$ such that if $w \neq 1$ in $G$, then
$\wtd{G}_w$ contains $G$ as a subgroup, while if $w=1$ then $\wtd{G}_w$ is
trivial. Taking the product $\wtd{G}_w \times G_+$ gives a group which either
contains $G$ as a subgroup, or is equal to $G_+$.  If we start with a group
$G_0 = \ang{S_0 : R_0}$ and word $w \in \mcF_{S_0}$, and apply the construction
to $G_0 \times G_-$, we get a group $H_w$ which is equal to $G_+$ when $w=1$,
or contains $G_-$ when $w \neq 1$. As a result, $H_w$ has (P) if and only if
$w$ is trivial in $G_0$. Taking $G_0$ to be a group with an undecidable word
problem yields the theorem.

The cocommutative compact quantum groups are exactly the group algebras $\C G$
for $G$ a discrete group.  As a result, the Adian--Rabin theorem for groups
implies an Adian--Rabin theorem for finitely presented cocommutative compact
quantum groups. We prove this in \Cref{sec:ccqg} for completeness.  However,
this doesn't immediately imply \Cref{thm:main}, because not every group algebra
is the quantum automorphism group of a finite graph
\cite{Banica-McCarthy}.
To prove \Cref{thm:main}, we use solution groups, which are the finite colimits
of finite abelian groups.  Solution groups are used in quantum information
theory to encode perfect strategies for linear system nonlocal games
\cite{Cleve-Liu-Slofstra}.  Undecidability results for solution groups have
been used to show that it's undecidable to determine if a linear system
nonlocal game has a perfect quantum strategy
\cite{Slofstra-Tsirelson,Slofstra-closed}, and subsequently to show that it's
undecidable to determine if two graphs are quantum isomorphic
\cite{Atserias-et-al}. Cai, Maran, and Young \cite{Cai-Maran-Young} use
\cite{Atserias-et-al} to show that it is undecidable to determine if there is a
``planar edge gadget'' whose signature differs at two fixed diagonal entries. 
\footnote{A quantum automorphism group $\Qut(X)$ is trivial if and only if
there is a planar edge gadget whose signature has distinct diagonal entries, so
this problem is closely related to deciding whether $\Qut(X)$ is trivial,
although not equivalent \emph{a priori}.}
Roberson and Schmidt \cite{Roberson-Schmidt} showed
that the group algebra of any solution group can be realized as the quantum
automorphism group of a vertex- and edge-coloured graph, and used this to
construct the first known graphs whose quantum automorphism groups are
non-classical but finite-dimensional, among other results.  More recently, the first
author, Roberson, and Schmidt \cite{Dobben-Roberson-Schmidt} proved that the
group algebra of any solution group can be realized as the quantum automorphism
group of an uncoloured graph.  We show that there is an analogue of the
Adian--Rabin construction for solution groups.  This construction is based on
the discovery of a simple solution group in \cite{Dobben-Roberson-Schmidt}.
Using this construction, we prove an Adian--Rabin theorem for solution groups,
which allows us to prove \Cref{thm:main}. This also allows us to show that it
is undecidable to determine if a solution group is trivial, finite, or abelian,
which was not known before.

Although it isn't quite a Markov property, our Adian--Rabin construction also
implies that it's undecidable to determine if a solution group has a
non-trivial finite-dimensional representation.  Translated over to quantum
automorphism groups, this implies:
\begin{theorem}\label{thm:main2}
    It is undecidable to determine, given a graph $X$ as input, whether
    $\Qut(X)$ has a non-trivial finite-dimensional $*$-representation.
\end{theorem}
In addition to the compact quantum automorphism group $\Qut(X)$, every finite
graph $X$ also has a discrete quantum automorphism group $\Qut_\delta(X)$,
which is built out of the finite-dimensional $*$-representations of $\Qut(X)$.
In particular, $\Qut_\delta(X)$ is trivial if and only if $\Qut(X)$ has no
non-trivial finite-dimensional $*$-representations, so \Cref{thm:main2} can be
restated as saying that it is undecidable to determine whether $\Qut_{\delta}(X)$
is trivial. This implies that there is a graph $X$ such that
$\Qut(X)$ is non-trivial but $\Qut_\delta(X)$ is trivial, answering a
question raised by Rollier and Vaes \cite[Remark~2.34]{Rollier-Vaes}.
We also provide a direct construction of such an example, without going
through a computability argument. Furthermore, we prove that it remains
undecidable to determine whether or not $\Qut(X)$ is trivial when we
are promised that $\Qut_\delta(X)$ is trivial.

\subsection*{Acknowledgements}
JvDdB is grateful to Jonas Eidesen, David Roberson, Lukas Rollier, Simon
Schmidt, and Christian Voigt for helpful discussions. Part of this paper was
written while the authors were in residence at Institut Mittag-Leffler in
Djursholm, Sweden during the research semester ``Operator Algebras and Quantum
Information'' in the spring of 2026, supported by the Swedish Research Council
under grant no.~2021-06594, and completed in the fall of 2026 while WS was
hosted at the Institut Henri Poincaré (IHP), supported by the CNRS, IHP (UAR
839 CNRS-Sorbonne Universite), and LabEx CARMIN (ANR-10-LABX-59-01). JvDdB was
supported by GA\v{C}R grant 25-17377S.  WS was supported by NSERC DG
2025-04395. 

\subsection*{AI disclosure}
The bulk of this paper was completed without AI assistance, and indeed the
key points were understood before the advent of the current generation of
reasoning models. We acknowledge the incidental use of AI models (ChatGPT and
Gemini) for e.g.~learning background material and searching the mathematical
literature. We report two occasions where AI use had a non-incidental impact
on the paper: we asked ChatGPT (Sol 5.6) whether there is a finitely presented
group algebra $\C G$ where $G$ is not finitely presentable, and (Astra) whether
the free symmetric group $\mathbb{S}_n^+$ ever contains cocommutative Hopf
$*$-subalgebras. In both cases, the answers turned out to be well-understood
in the literature, but we may have overlooked them or taken longer to find them
if we did not have AI assistance. 

AI models were not used for writing assistance.

\section{Solution groups}
\label{sec:solgp}

In this section, we prove our Adian--Rabin theorem for solution groups.
We start with the definition of a (homogeneous) solution group.

\begin{defn}[\cite{Roberson-Schmidt}]\label{def:solngrp}
	Let $Z = \{z_1,\ldots,z_n\}$ be an ordered set of size $n \geq 1$, and $A$ be an
	$m \times n$ matrix over $\Z_2$ for some $m \geq 0$. The \emph{(homogeneous)
	solution group of $A$} with generating set $Z$ is the finitely
	presented group
	\begin{align*}
		\Gamma(A,Z) = \ang{\ Z\ :\ &  z_j^2=1 \text{ for all } 1 \leq j \leq n, \\
		                           &  \prod_{j=1}^n z_j^{A_{ij}} = 1 \text{ for all } 1 \leq i \leq m, \\
		                           &  [z_j,z_k]=1 \text{ if } A_{ij}, A_{ik} \neq 0 \text{ for some }1 \leq i \leq m}.
	\end{align*}
	We let $R(A,Z) \subseteq \mcF_Z$ denote the set of defining relations of $\Gamma(A,Z)$.
	We set
	\begin{equation*}
		\Gamma(A) := \Gamma(A,\{x_1,\dots,x_n\}) \text{ and } R(A) := R(A,\{x_1,\dots,x_n\}).
	\end{equation*}
\end{defn}
Here $[a,b] := a b a^{-1} b^{-1}$ is the group commutator, and $\mcF_Z$ is the free group generated by $Z$.
Note that $\Gamma(A,Z)$ is defined as a presented group, meaning that the generating set
$Z$ and relations $R(A,Z)$ are part of the data of the group.  The notation
$\Gamma(A,Z)$ and $R(A,Z)$ is chosen to make it convenient to use the matrix
$A$ with different generating sets; we use $\{x_1,\ldots,x_n\}$ as the default
set of generators. Note from this definition that $\Gamma(A_0 \oplus A_1)$ is
the free product $\Gamma(A_0) * \Gamma(A_1)$.

There is also a non-homogeneous solution group defined as follows:
\begin{defn}[\cite{Cleve-Liu-Slofstra}]
	Let $Z = \{z_1,\ldots,z_n\}$ be an ordered set of size $n \geq 1$, and $Ax=b$ be an
	$m \times n$ linear system over $\Z_2$ for some $m \geq 0$. The \emph{solution group
	of $Ax=b$} with generating set $Z$ is the finitely presented group
	\begin{align*}
		\Gamma(A,b,Z) = \ang{\ Z \cup \{J\} \ :\ &  z_j^2 = [z_j,J] = J^2 = 1 \text{ for all } 1 \leq j \leq n, \\
		                                         &  \prod_{j=1}^n z_j^{A_{ij}} = J^{b_i} \text{ for all } 1 \leq i \leq m, \\
		                                         &  [z_j,z_k]=1 \text{ if } A_{ij}, A_{ik} \neq 0 \text{ for some }1 \leq i \leq m}.
	\end{align*}
	We let $\Gamma(A,b) := \Gamma(A,b,\{x_1,\dots,x_n\})$.
\end{defn}
In this paper, solution group will mean a homogeneous solution group by default. Any
non-homogeneous solution group can be realized as a homogeneous solution group by
regarding $J$ as one of the variables, and adding dummy variables and relations to
force $J$ to commute with the other variables
(see e.g.~\cite[Lemma~5.4]{Dobben-Roberson-Schmidt} for more details).

We now show that there is an Adian--Rabin construction for solution groups.
For the theorem statement, let $\mcN(x)$ denote the normal subgroup generated by $x$.
\begin{thm}
	\label{thm:adianrabin}
	For every $m \times n$ matrix $A$ over $\Z_2$, there is an $\tilde{m} \times \tilde{n}$ matrix
	$\wtd{A}$ over $\Z_2$, where $\tilde{m} \geq m$ and $\tilde{n} \geq n$, such that
	\begin{enumerate}[(1)]
		\item there is a homomorphism $\iota : \Gamma(A) \to \Gamma(\wtd{A})$ sending $x_i \mapsto x_i$ for $1 \leq i \leq n$,
		
		\item if $x_1$ is non-trivial in $\Gamma(A)$, then the homomorphism $\iota$ is injective,
		
		\item\label{itm:diagram} if $\varphi : \Gamma(A) \to G$ is a homomorphism to a finite group
		$G$ such that $\varphi(x_1) \neq 1$, then there is a finite group $H$
		containing $G$ as a subgroup, and a homomorphism $\wtd{\varphi} :
		\Gamma(\wtd{A}) \to H$, such that the following diagram commutes:
		\begin{tikzcd}
			\Gamma(A) \arrow[r,"\varphi"] \arrow[d,"\iota"] & G \arrow[d,hook] \\
			\Gamma(\wtd{A}) \arrow[r,"\wtd{\varphi}"] & H
		\end{tikzcd}, and
		
		\item the quotient $\Gamma(\wtd{A}) / \mcN(x_1)$ is trivial.
	\end{enumerate}
	Furthermore, the mapping $A \mapsto \wtd{A}$ is computable in polynomial time.
\end{thm}

To prove \Cref{thm:adianrabin}, we use two specific solution groups. For the first, let
$A_\MP \cdot x = 0$ be the $\Z_2$-linear system on variables $x_1,\ldots,x_{10}$ which states
the sum of the first three columns and every row in the following table is zero:
\begin{center}
	\begin{tabular}{|c|c|c|c|}
		\hline
		$x_1$ & $x_2$ & $x_3$ & $x_{10}$ \\
		\hline
		$x_4$ & $x_5$ & $x_6$ & $x_{10}$ \\
		\hline
		$x_7$ & $x_8$ & $x_9$ & $x_{10}$ \\
		\hline
	\end{tabular}
\end{center}
Note that $x_{10}$ commutes with $x_1,\ldots,x_9$ in $\Gamma(A_\MP)$. If we rename $x_{10}$ to $J$, then
$\Gamma(A_\MP,\{x_1,\ldots,x_9,J\})$ is the homogeneous solution group $\Gamma(A,b)$,
where $Ax=b$ is the $\Z_2$-linear system $Ax=b$ on variables $x_1,\ldots,x_9$ which states
that every row in the following table sums to $1$, and every column to $0$:
\begin{align*}
	\begin{array}{|c|c|c|c|}
		\hline
		x_1 & x_2 & x_3 \\
		\hline
		x_4 & x_5 & x_6 \\
		\hline
		x_7 & x_8 & x_9 \\
		\hline
	\end{array}
	\tag*{$(*)$}
\end{align*}
This linear system is the famous Mermin--Peres magic square \cite{Peres,Mermin},
and it is well-known that $\Gamma(A_\MP,\{x_1,\ldots,x_9,J\})$
is finite, that all the generators $x_1,\ldots,x_9,J$ are non-trivial, and that
\begin{equation*}
	[x_i,x_j] = \begin{cases}
		1 & x_i \text{ and } x_j \text{ are in the same row or column of $(*)$} \\
		J & \text{otherwise} \\
	\end{cases}
\end{equation*}
for all $1 \leq i,j \leq 9$. In particular, $\Gamma(A_\MP) / \mcN(J)$ is abelian.

The other solution group we need is the finite simple solution group discovered in \cite{Dobben-Roberson-Schmidt}.
\begin{lemma}
	\label{lem:perfectsolngrp}
	There is a $1645 \times 315$ matrix $A_\simple$ over $\Z_2$, such that
	$\Gamma(A_\simple)$ is finite, non-trivial, and simple, and $|[x_1,x_2]| = 2$ in $\Gamma(A_\simple)$.
\end{lemma}
\begin{proof}
	Consider the matrix $M_{A_8}$ constructed from the alternating group $A_8$ using the construction from \cite[Definition~6.1]{Dobben-Roberson-Schmidt}.
	The columns of this matrix are indexed by the elements of order $2$ in $A_8$, and we know from \cite[Remark~6.9]{Dobben-Roberson-Schmidt} that $\Gamma(M_{A_8}) \cong A_8$ via the isomorphism $x_g \mapsto g$.
	Choose two elements $g_1,g_2 \in A_8$ of order $2$ such that $|[g_1,g_2]| = 2$.
	For example, we may take $g_1 = (12)(34)$ and $g_2 = (35)(46)$, so that $[g_1,g_2] = (34)(56)$.
	Reordering the columns so that the first two are indexed by $g_1$ and $g_2$ gives the desired conclusion.
\end{proof}

We now describe the main construction for the proof of \Cref{thm:adianrabin}.
For the construction, we say that a finitely presented group $P = \langle S : R \rangle$ is \emph{generated by involutions} if $s^2 = 1$ in $P$ for all $s \in S$.

\begin{construction}
	Let $P = \langle S : R \rangle$ be a finitely presented group generated by involutions, where $S = \{s_1,\ldots,s_n\}$.
	Define a finitely presented group $\AR(P) = \langle \AR(S) : \AR(S,R) \rangle$ by
	\begin{alignat*}{2}
		\AR(S) &= \mathrlap{S \cup S' \cup \bigcup_{j \in C_0} Y_j \cup \bigcup_{(j,k) \in C_0 \times C_1} Z_{jk} \quad \text{and}} \\[1ex]
		\AR(P) &= \ang{\AR(S) \ :\ && R \cup \bigcup_{j \in C_0}  R(A_\MP,Y_j) \cup \bigcup_{(j,k) \in C_0 \times C_1} R(A_\simple, Z_{jk}), \\
		                            &&& s_{i} s_{i0} s_{i1} = s_{i0}^2 = s_{i1}^2 = [s_{i0},s_{i1}]=1 \text{ for all } 2 \leq i \leq n,\\
		                            &&& s_j = y_{j1} \text{ for all } j \in C_0, \\
		                            &&& y_{jk_1} = z_{j(k_1,k_2)1} \text{ and } y_{jk_2} = z_{j(k_1,k_2)2} \text{ for all } (j,(k_1,k_2)) \in C_0 \times C_1 },
	\end{alignat*}
	where $A_\MP$ and $A_\simple$ are the two matrices mentioned above, $r$ is the number of columns of $A_\simple$, $C_0 = \{1\} \cup \{(j_1,j_2) : 2 \leq j_1 \leq n, j_2 \in \{0,1\}\}$, $C_1$ is the set of pairs $(k_1,k_2)$ with $1 \leq k_1,k_2 \leq 9$ such that $x_{k_1}$ and $x_{k_2}$ are not in the same row or column in $(*)$, $S' = \{s_{j} : j \in C_0\}$, $Y_j = \{y_{j1},y_{j2},\ldots,y_{j9}, J\}$ for $j \in C_0$,
	and $Z_{jk} = \{z_{jk\ell}\}_{\ell=1}^r$ for $(j,k) \in C_0 \times C_1$.
\end{construction}

\begin{lemma}
	\label{lem:inclusionhom}
	For every finitely presented group $P$ generated by involutions, there is
	a homomorphism $\iota_P : P \to \AR(P)$ sending $s \mapsto s$ for all
	$s \in S$.  Moreover, if $P = \langle S : R \rangle$ and $P' =
	\langle S : R' \rangle$ are finitely presented groups generated by
	involutions on the same generating set $S$, and $\varphi : P \to P'$ is
	a homomorphism sending $s \mapsto s$ for all $s \in S$,
	then there is a homomorphism $\AR(\varphi) : \AR(P) \to \AR(P')$
	sending $s \mapsto s$ for all $s \in \AR(P)$, and this homomorphism
	satisfies $\iota_{P'} \varphi = \AR(\varphi) \iota_P$.
\end{lemma}
\begin{proof}
	Since the presentation of $\AR(P)$ contains the relations of $P$, the
	mapping $s \mapsto s$ for all $s \in S$ induces a homomorphism $P \to \AR(P)$.  Since
	$P$ and $P'$ have the same generating set $S$, $\AR(P)$ and $\AR(P')$
	have the same generating set $\AR(S)$. Also, with the exception of the 
	relations from $R$ and $R'$, $\AR(P)$ and $\AR(P')$ have the same defining
	relations, i.e. $\AR(S,R') = (\AR(S,R) \setminus R) \cup R'$. Since $\varphi$
	is a homomorphism, we conclude that there is a homomorphism $\AR(\varphi) :
	\AR(P) \to \AR(P')$ sending $s \mapsto s$ for all $s \in \AR(S)$, and
	this homomorphism necessarily satisfies $\iota_{P'} \varphi = \AR(\varphi) \iota_{P}$.
\end{proof}

\begin{lemma}
	\label{lem:amalgamated}
	If $P$ is a finitely presented group generated by involutions $S = \{s_1,\ldots,s_n\}$ such that $s_1 \neq 1$ in $P$,
	then there is a sequence of finite groups $K_1,\ldots,K_t$ and subgroups $H_i \subseteq K_i$ such that $\AR(P)$ is an iterated
	amalgamated free product
	\[ \AR(P) = \left( \left( ((P \times \Z_2) *_{H_1} K_1) *_{H_2} K_2\right) \cdots \right) *_{H_t} K_t. \]
	Consequently, the map $\iota_P$ is injective whenever $s_1 \neq 1$.
\end{lemma}
\begin{proof}
	Suppose $P = \ang{ S : R}$. Let 
	\begin{align*}
		\AR_0 = \ang{S \cup \{J\} : R, [s_i,J] = J^2 = 1 \text{ for all } s_i \in S}
	\end{align*}
	and
	\begin{align*}
		\AR_1 = \ang{S \cup \{J\} \cup S' \ :\ & R, [s_j,J] = J^2 = 1 \text{ for all } j \in C_0, \\
		& s_{i} s_{i0} s_{i1} = s_{i0}^2 = s_{i1}^2 = [s_{i0},s_{i1}]=1 \text{ for all } 2 \leq i \leq n}.
	\end{align*}
	Clearly $\AR_0 = P \times \Z_2$. Since $s_i = s_{i0} s_{i1}$ in $\AR_1$ for all $2 \leq i \leq n$,
	we have $[J,s_i] = 1$ in $\AR_1$ for all $1 \leq i \leq n$. If $s_i \neq 1$ in
	$P$ for some $2 \leq i \leq n$, then adding $s_{i0}$ and $s_{i1}$ to the
	presentation of $\AR_1$ is equivalent to taking the amalgamated free
	product with $\ang{s_{i0}, s_{i1}, J : [s_{i0},s_{i1}] = [s_{i0},J] =
	[s_{i1},J] = s_{i0}^2 = s_{i1}^2 = J^2 = 1} \iso \Z_2^3$ over the subgroup
	$\ang{s_i , J} = \ang{s_{i_0} s_{i_1}, J} \iso \Z_2 \times \Z_2$. If $s_i=1$,
	then adding $s_{i0}$ and $s_{i1}$ is equivalent to taking the amalgamated
	free product with $\ang{s_{i0}, s_{i1}, J : s_{i0} s_{i1} = [s_{i0}, J] = s_{i0}^2 = J^2 = 1 } \iso \Z_2 \times \Z_2$
	over the subgroup $\ang{J} \iso \Z_2$. Thus $\AR_1$ is an iterated amalgamated
	free product $(P \times \Z_2) *_{H_1} K_1 *_{H_2} \cdots *_{H_{n-1}} K_{n-1}$, where
	$(K_{i},H_i) = (\Z_2^3,\Z_2^2)$ if $s_{i+1} \neq 1$ in $P$, and $(K_i,H_i) = (\Z_2^2,\Z_2)$ otherwise.
	Note that $s_{ij} \neq 1$ in $K_{i-1}$, and therefore $s_{ij} \neq 1$ in $\AR_1$ for all $2 \leq i \leq n$, $j \in \{0,1\}$.
	
	Next, we look at the group
	\begin{align*}
		\AR_{2} = \ang{S \cup S' \cup \bigcup_{j \in C_0} Y_j \ :\ & R \cup \bigcup_{j \in C_0} R(A_\MP,Y_j), \\
		& s_{i} s_{i0} s_{i1} = s_{i0}^2 = s_{i1}^2 = [s_{i0},s_{i1}]=1 \text{ for all } 2 \leq i \leq n, \\
		& s_j = y_{j1} \text{ for all } j \in C_0}.
	\end{align*}
	Observe that the set $\{y_{j\ell} : j \in C_0, 1 \leq \ell \leq 9\}$ generates $\AR_2$.
	Since $[J,y_{j\ell}]=1$ for all $j \in C_0$ and $\ell=1,\ldots,9$, we see that $J$ continues to be central in $\AR_2$.
	Since $s_1 \neq 1$ in $P$ by hypothesis, $\ang{s_j,J}_{\AR_1} \iso \Z_2 \times \Z_2
	\iso \ang{y_{j1},J}_{\Gamma(A_\MP,Y_j)}$ for all $j \in C_0$, so
	$\AR_2$ is the iterated amalgamated free product of $\AR_1$ with a copy of $\Gamma(A_\MP,Y_j)$ for
	every $j \in C_0$. In particular, this implies that $\ang{Y_{j}}_{\AR_2}
	\iso \Gamma(A_\MP,Y_j)$. Since $|[x_1,x_2]|=2$ in $\Gamma(A_\simple)$,
	we have $\ang{x_1,x_{2}}_{\Gamma(A_\simple)} \iso D_4$, the dihedral group of order $8$.
	The same is true of $\ang{x_{k_1},x_{k_2}}_{\Gamma(A_\MP)}$ for $(k_1,k_2) \in C_1$.
	Hence $\AR(P)$ is an iterated amalgamated free product of $\AR_2$ with a
	copy of $\Gamma(A_\simple,Z_{jk})$ for each $(j,k) \in C_0 \times C_1$.
\end{proof}

\begin{lemma}
	\label{lem:trivialquotient}
	If $P$ is a finitely presented group generated by involutions $S = \{s_1,\ldots,s_n\}$ such that $s_1 = 1$ in $P$,
	then $\AR(P)$ is trivial.
\end{lemma}
\begin{proof}
	If $s_1 =1$ in $P$, then $s_1 = 1$ in $\AR(P)$. This means that $y_{11} = 1$, so the
	relations of $R(A_\MP,Y_{1})$ give that $J = [y_{11},y_{15}] = 1$. Hence the relations in
	$R(A_\MP,Y_j)$ imply that $[y_{jk_1},y_{jk_2}]=1$ for all $(j,k) \in C_0
	\times C_1$, where $k = (k_1,k_2)$. But this means that $[z_{jk1},z_{jk2}]=1$.  Since
	$\Gamma(A_\simple,Z_{jk})$ is simple and $[z_{jk1},z_{jk2}] \neq 1$ in
	$\Gamma(A_\simple,Z_{jk})$, setting $[z_{jk1},z_{jk2}]=1$ implies that
	$z_{jk\ell}=1$ for all $\ell=1,\ldots,r$.  Tracing backwards, we get that
	$y_{j\ell}=1$ for $\ell=1,\ldots,9$, that $s_j = 1$ for all $j \in C_0$, and finally that $s_i=1$ for $2 \leq i
	\leq n$. Since all the generators are identity, we conclude that $\AR(P)$
	is trivial.
\end{proof}

\begin{lemma}
	\label{lem:itslinear}
	Let $A$ be an $m \times n$ matrix over $\Z_2$, and let $S =
	\{s_1,\ldots,s_n\}$.  Then there is an $\tilde{m} \times \tilde{n}$ matrix
	$\wtd{A}$ over $\Z_2$, where $\tilde{m} \geq m$ and $\tilde{n} \geq n$,
	such that $\AR(\Gamma(A,S)) = \Gamma(\wtd{A},\AR(S))$. The matrix $\wtd{A}$
	is constructible from $A$ in polynomial time. 
\end{lemma}
\begin{proof}
	Since $s_{i0}$ and $s_{i1}$ commute and $s_i$ is their product, we have that $s_i$ commutes with $s_{i0}$ and $s_{i1}$ in $\AR(\Gamma(A,S))$.\footnote{Actually, the relation $[s_{i0},s_{i1}] = 1$ is superfluous, since the relations $s_i^2 = s_{i0}^2 = s_{i1}^2 = 1$
	and $s_i s_{i0} s_{i1} = 1$ imply that $s_i$, $s_{i0}$ and $s_{i1}$ commute, but we keep this relation in the presentation for clarity.}
	Thus we can recognize
	all the relations of $\AR(\Gamma(A,S))$ as linear relations as in \Cref{def:solngrp},
	and find a matrix $\wtd{A}$ such that $\AR(\Gamma(A,S)) = \Gamma(\wtd{A},\AR(S))$.
	Since the number of new relations in $\AR(\Gamma(A,S))$ is linear in $n$,
	the matrix $\wtd{A}$ can be constructed from $A$ in polynomial time.
\end{proof}

\begin{proof}[Proof of \Cref{thm:adianrabin}]
	Suppose $A$ is an $m \times n$ matrix over $\Z_2$, and let $X =
	\{x_1,\ldots,x_n\}$.  By \Cref{lem:itslinear}, we can construct a matrix
	$\wtd{A}$ (in polynomial time) such that $\Gamma(\wtd{A},\AR(X)) =
	\AR(\Gamma(A))$. Parts (1) and (2) of the theorem follow from Lemmas
	\ref{lem:inclusionhom} and \ref{lem:amalgamated}.
	For part (4), consider $\Gamma(A) / \mcN(x_1)$ presented as $\langle X : R(A,X) \cup \{x_1\} \rangle$.
	By construction, we have $\AR(\Gamma(A) / \mcN(x_1)) = \langle \AR(X) : \AR(X, R(A,X) \cup \{x_1\}) \rangle = \langle \AR(X) : \AR(X,R(A,X)) \cup \{x_1\} \rangle = \AR(\Gamma(A)) / \mcN(x_1)$.
	Since $x_1 = 1$ in $\AR(\Gamma(A) / \mcN(x_1))$, $\AR(\Gamma(A)) / \mcN(x_1)$ is trivial by \Cref{lem:trivialquotient}.
	
	This leaves part (3). Suppose we are given a homomorphism $\varphi : \Gamma(A)
	\to G$ where $G$ is finite and $\varphi(x_1) \neq 1$.  Let $G_0 = \varphi(\Gamma(A))$ be the image of
	$\varphi$, and take a presentation {${G_0 = \ang{S_0 : R_0}}$} of $G_0$ with
	generating set $S_0 = \{t_1,\ldots,t_n\}$, where $t_i = \varphi(x_i)$ for $1 \leq i \leq n$.
	Evidently $G_0$ is generated by involutions. By \Cref{lem:inclusionhom}, there is a homomorphism
	$\AR(\varphi) : \AR(\Gamma(A)) \to \AR(G_0)$ such that $\AR(\varphi) \iota_{\Gamma(A)}
	= \iota_{G_0} \varphi$.
	Since $\varphi(x_1) \neq 1$, it follows from \Cref{lem:amalgamated} that $G_0$ is a subgroup of $\AR(G_0)$.
	Set $H' = G *_{G_0} \AR(G_0)$, and
	let $\psi'$ be the composition of $\AR(\varphi)$ with the inclusion $\AR(G_0) \incl H'$.
	Then $H'$ contains $G$ as a subgroup, and the diagram
	\[ \begin{tikzcd}
		\Gamma(A) \arrow[r,"\varphi"] \arrow[d,"\iota"] & G \arrow[d,hook] \\
		\AR(\Gamma(A)) \arrow[r,"\psi'"] & H'
	\end{tikzcd} \]
	commutes.
	By \Cref{lem:amalgamated}, $\AR(G_0)$ is an iterated amalgamated free product of $G_0 \times \Z_2$ with a
	sequence of finite groups, and hence the same is true of $H'$. By a result
	of Baumslag \cite[Theorem~2]{Baumslag-RF}, $H'$ is residually finite. Since
	$G$ is finite, there is a finite group $H$ and a homomorphism $\psi'' : H' \to H$
	such that $\psi''|_G$ is injective, meaning that $H$ contains $G$ as a
	subgroup. Setting $\psi = \psi'' \circ \psi'$ gives the desired
	homomorphism $\Gamma(\wtd{A},\AR(X)) = \AR(\Gamma(A)) \to H$.
\end{proof}

\Cref{thm:adianrabin} leads to the following analogue of the Adian--Rabin theorem
for solution groups.
\begin{corollary}
	\label{cor:markovprop}
	Let $(P)$ be a Markov property of groups, such that there is a solution
	group $\Gamma_+$ with property $(P)$. Then it is undecidable (specifically, RE-hard) to determine
	if a solution group has $(P)$.
\end{corollary}
\begin{proof}
	Let $\Gamma(A_+)$ be the solution group with property $(P)$, let $G_-$ be
	a finitely presented group such that any group containing $G_-$ as a subgroup does not have
	property $(P)$, and let $M$ be a Turing machine with an undecidable
	halting problem. By \cite[Theorem~3.1]{Slofstra-Tsirelson}, there is a solution group
	$\Gamma(A_-)$ containing $G_-$ as a subgroup. By \cite[Corollary~3.3]{Slofstra-Tsirelson},
	there is a computable function from bit strings $w \in
	\{0,1\}^*$ to $\Z_2$-linear systems $T_w x = b_w$, such that $J=1$ in
	$\Gamma(T_w,b_w)$ if and only if $M$ halts on input $w$. Since every
	non-homogeneous solution group can be turned into a homogeneous solution
	group, we can turn this into a computable function from bitstrings $w$ to
	matrices $A_w$ such that $x_1 = 1$ in $\Gamma(A_w)$ if and only if $M$
	halts on input $w$. Since $\Gamma(A_w \oplus A_-) = \Gamma(A_w) * \Gamma(A_-)$,
	we have $x_1 =1$ in $\Gamma(A_w \oplus A_-)$ if and only if $M$ halts on $w$.
	
	Let $\wtd{A}_w$ be the matrix constructed from $A = A_w \oplus A_-$ in
	\Cref{thm:adianrabin}, and consider the family of matrices $\{ \wtd{A}_w \oplus A_+ \}_{w \in \{0,1\}^*}$.
	If $M$ halts on input $w$, then by part (a) of the
	theorem, $x_1 = 1$ in $\Gamma(\wtd{A}_w)$, and by part (d),
	$\Gamma(\wtd{A}_w)$ is trivial. Hence $\Gamma(\wtd{A}_w \oplus A_+)
	= \Gamma(A_+)$ has property $(P)$. If $M$ does not halt on input $w$, then
	$x_1 \neq 1$ in $\Gamma(A_w \oplus A_-)$, so $\Gamma(\wtd{A}_w)$ contains
	$\Gamma(A_w \oplus A_-)$ as a subgroup. Hence $\Gamma(\wtd{A}_w \oplus A_+)$ contains
	$G_0$ as a subgroup, and does not have property $(P)$. Thus, if we could decide
	whether a solution group has property $(P)$, then we could decide the halting
	problem for $M$.
\end{proof}

\begin{example}
	\label{xmpl:Markov}
	Examples of Markov properties include being trivial, finite, and abelian. For each
	of these properties, taking $A$ to be the $0 \times 0$ empty matrix, we get a solution group $\Gamma(A)$ which is trivial, and hence finite and abelian.
	Thus \Cref{cor:markovprop} applies to these properties.
\end{example}

The fact that it is undecidable to determine whether a solution group is trivial
implies that there is an infinite solution group which does
not have a non-trivial finite-dimensional representation. Indeed, if every
non-trivial solution group had a non-trivial finite dimensional representation,
then it would be possible to determine whether a solution group $\Gamma$ is trivial
by searching simultaneously for a proof that $\Gamma$ is trivial, and for a non-trivial
finite-dimensional representation of $\Gamma$. Thus there must be a non-trivial
solution group with no non-trivial finite-dimensional representations, and such a
group must be infinite. Using \cref{thm:adianrabin}, we can also give a direct
construction of such an example:
\begin{cor}\label{cor:nofindim}
    There is an infinite solution group with no non-trivial finite-dimensional
    representations.
\end{cor}
\begin{proof}
        By \cite[Corollary~3.3]{Slofstra-Tsirelson}, there is a non-homogeneous
        $\Z_2$-linear system $T x = b$ such that $J \neq 1$ in $\Gamma(T,b)$, but
        $J$ is trivial in all finite-dimensional representations of
        $\Gamma(T,b)$. As in the proof of \Cref{cor:markovprop}, this means
        that there is a matrix $A$ such that $x_1 \neq 1$ in $\Gamma(A)$, but
        $x_1$ is trivial in all finite-dimensional representations.
        Applying \Cref{thm:adianrabin}, we get a matrix $\wtd{A}$ such that
        $\Gamma(\wtd{A})$ contains $\Gamma(A)$ as a subgroup, but 
        $\Gamma(\wtd{A}) / \mcN(x_1)$ is trivial. We conclude that
        $\Gamma(\wtd{A})$ is non-trivial, but has no non-trivial
        finite-dimensional representations (and in particular, is infinite).
\end{proof}

Let (P) be the property of having a non-trivial finite-dimensional representation.
Neither (P) nor its complement is a Markov property, since any finitely presented group $G$ is contained in a finitely presented group which has a non-trivial finite-dimensional representation (e.g. $G \times \Z_2)$, and also in a finitely presented group which has no non-trivial finite-dimensional representations.\footnote{One way to show this is with the Adian--Rabin construction.
Suppose we're given a group $G$.
Take a finitely presented group $H$ with an element $x \neq 1$ which is trivial in every finite-dimensional representation.
By the Adian--Rabin construction, there is a finitely presented group $K$ containing $G \times H$ such that $K / \mcN(x)$ is trivial.
Then $K$ has no non-trivial finite-dimensional representations.}
Nevertheless, using our Adian--Rabin construction and a different undecidability result for solution groups from \cite{Slofstra-closed}, we can show that this property is also undecidable.

\begin{corollary}
	\label{cor:findimar}
	It is undecidable to determine if a solution group has a non-trivial finite-dimensional representation.
\end{corollary}
\begin{proof}
	A theorem of Mal'cev \cite{Malcev} states that for any finitely generated group $G$ and element $x \in G$, there is a
	finite-dimensional unitary representation $\varphi : G \to U(\C^d)$ for some $d$ with
	$\varphi(x) \neq 1$ if and only if there is a homomorphism $\psi : G \to H$ for some finite
	group $H$ such that $\psi(x) \neq 1$. In particular, having a non-trivial finite-dimensional
	representation is the same as having a non-trivial homomorphism to a finite group.
	
	It is shown in \cite{Slofstra-closed} that it is undecidable to determine, given a non-homogeneous solution
	group $\Gamma(T,b)$, whether $J$ is non-trivial in some finite-dimensional representation.
	As mentioned previously, every non-homogeneous solution group $\Gamma(T,b)$ can be realized as a homogeneous solution group,
	meaning that there is a matrix $A$ and an isomorphism $\Gamma(T,b) \to \Gamma(A)$ sending $J \mapsto x_1$.
	Since the map $(T,b) \mapsto A$ is computable, it is undecidable to determine, given a solution group $\Gamma(A)$, whether $x_1$ is non-trivial in a finite-dimensional of $\Gamma(A)$.
	Equivalently (by Mal'cev's theorem), it is undecidable to determine whether there is a finite group $H$ and a homomorphism $\varphi : \Gamma(A) \to H$ such that $\varphi(x_1) \neq 1$.
	
	Suppose $\Gamma(A)$ is homogeneous solution group, and let $\wtd{A}$ be the matrix from \Cref{thm:adianrabin}.
	If there is a finite group $G$ and a homomorphism $\varphi : \Gamma(A) \to G$ such that $\varphi(x_1) \neq 1$,
	then by part (c) of \Cref{thm:adianrabin}, there is a finite group $H \supseteq G$ and a homomorphism $\psi : \Gamma(\wtd{A}) \to H$
	such that $\psi(x_1) = \varphi(x_1) \neq 1$. Hence $\Gamma(\wtd{A})$ has a
	non-trivial homomorphism to a finite group.
	Conversely, suppose that $\psi : \Gamma(\wtd{A}) \to H$ is a homomorphism to a finite group $H$.
	If $\psi(x_1) = 1$, then $\psi$ factors through $\Gamma(\wtd{A}) / \mcN(x_1) = \{1\}$, so
	$\psi$ is trivial. Thus, if $\Gamma(\wtd{A})$ has a non-trivial homomorphism $\psi$ to a finite
	group, then $\psi(x_1) \neq 1$, in which case $\varphi = \psi|_{\Gamma(A)} : \Gamma(A) \to H$ is also
	a homomorphism to a finite group with $\varphi(x_1) \neq 1$. We conclude that
	there is a homomorphism $\varphi : \Gamma(A) \to G$ to a finite group $G$ such that $\varphi(x_1) \neq 1$
	if and only if $\Gamma(\wtd{A})$ has a non-trivial homomorphism to a finite group.
	
	Since the mapping $A \mapsto \wtd{A}$ is computable, it is undecidable to
	determine if a solution group has a non-trivial finite-dimensional representation.
\end{proof}

Note that the property (P) of having a non-trivial finite-dimensional
representation is in RE.  However, the proof in \cite{Slofstra-closed} that
it's undecidable to determine whether $J$ is non-trivial in some
finite-dimensional representation is not from a reduction from the halting
problem, so we do not know if (P) is RE-complete.

\section{Cocommutative compact quantum groups}\label{sec:ccqg}

Our goal is to prove an Adian--Rabin theorem for quantum automorphism groups of
graphs. In this section, we recall some background about compact quantum
groups, including that a cocommutative compact quantum group is a group
algebra. Although we don't need it for the next section, we observe that the
usual Adian--Rabin theorem for groups lifts to cocommutative compact quantum
groups.

A \emph{$*$\nobreakdash-algebra} is an algebra over $\C$ with a conjugate linear
map $* : A \to A$ such that $(ab)^* = b^*a^*$ and $(a^*)^* = a$ for all
$a, b \in A$. All the algebras we consider will be unital. If $A$ and $B$ are
$*$-algebras, then $A \tensor B$ is a $*$-algebra with $(a \tensor b)^* = a^*
\tensor b^*$. Given a set $X$, we let $\C^*\ang{X}$ be the $*$-algebra of
noncommutative $*$-polynomials with indeterminates $X$. If $R \subseteq
\C^*\ang{X}$, we let ${\C^*\ang{X : R}} = \C^*\ang{X} / \ang{R}$, where
$\ang{R}$ is the two-sided $*$-ideal generated by $R$.  A \emph{finitely
presented $*$-algebra} is a $*$-algebra $A$ of the form $A = \C^*\ang{X : R}$
for finite sets $X$ and $R$, where $X$ and $R$ are remembered as part of the
data of $A$. If $\K$ is a subfield of $R$ which is closed under conjugation, we
let $\K^*\ang{X} \subseteq \C^*\ang{X}$ denote the $*$-subring of $*$-polynomials
with coefficients in $\K$.  We say that $\C^*\ang{X:R}$ is \emph{presented over $\K$}
if $R \subseteq \K^*\ang{X}$.  If $A = \C^*\ang{X : R}$ and $B = \K^*\ang{Y :
S}$ are finitely presented, then the tensor product $A \tensor B$ is finitely
presented with presentation $\C^*\ang{X \sqcup Y , R \sqcup S \sqcup \{xy - yx
: x \in X , y \in Y\}}$. A $*$-algebra is \emph{finitely presentable (over $\K$)} if it
has a finite presentation (over $\K$). A \emph{finitely presented $*$-homomorphism} $A \to
B$ is a pair $(\varphi , \wtd\varphi)$ where $\varphi : A \to B$ is a
$*$-homomorphism and $\wtd\varphi : X \to \C^*\ang{Y}$ is a function such that
$\varphi(x) = \wtd\varphi(x) + \langle S \rangle$ for all $x \in X$. Usually we
drop $\wtd\varphi$ from the notation, and use the convention that $\wtd{\varphi}$
is remembered as part of the data of $\varphi$. If $\wtd\varphi(X) \subseteq
\K^*\ang{Y}$, then we say that $\varphi$ is \emph{finitely presented over $\K$}.

A \emph{Hopf algebra} is a (unital) algebra $A$ together with an associative
coproduct $\Delta : A \to A \tensor A$, a counit $\epsilon : A \to \C$, and an
antipode $S : A \to A$, which satisfy some compatibility conditions.  A
\emph{Hopf $*$-algebra} is a Hopf algebra where $A$ is additionally a
$*$-algebra and $\Delta$ and $\epsilon$ are $*$-homomorphisms.
A \emph{Hopf $*$-subalgebra} of a Hopf $*$-algebra $A$ is a $*$-subalgebra $B$
such that $\Delta(B) \subseteq B \otimes B$ and $S(B) \subseteq B$.
A \emph{compact quantum group} is a Hopf $*$-algebra which admits a \emph{Haar
state}, meaning a linear map $h : A \to \C$ such that $h(1) = 1$ and $(\id_A
\tensor h)\Delta(a) = h(a)\cdot 1_A = (h \tensor \id_A)\Delta(a)$ for all $a
\in A$.  The Haar state is positive and faithful, meaning that $h(a^*a) > 0$
for all $a \in A \setminus \{0\}$. If $A$ and $B$ are compact quantum groups
with coproducts $\Delta_A$ and $\Delta_B$ respectively, then the free product
algebra $A \Ast B$ can be made into a compact quantum group by taking the
coproduct to be the unique $*$-homomorphism $\Delta : A \Ast B \to (A \Ast B)
\otimes (A \Ast B)$ such that $\Delta \circ i_A = (i_A \otimes i_A) \circ
\Delta_A$ and $\Delta \circ i_B = (i_B \otimes i_B) \circ \Delta_B$, where
$i_A$ and $i_B$ are the inclusions of $A$ and $B$ inside $A \Ast B$.  We refer
to \cite{Timmermann} for more background on Hopf $*$-algebras.

We'll say that a Hopf $*$-algebra $A$ is \emph{finitely presented (over $\K
\subseteq \C$)} if $A$ is a finitely presented $*$-algebra (over $\K$) and
$\Delta$ is a finitely presented $*$-homomorphism (over $\K$). The counit
and antipode of a Hopf $*$-algebra are uniquely determined by the coproduct,
so we don't ask for these as part of the presentation.  If $A = \C^*\ang{X :
R}$ and $B = \C^*\ang{X' : R'}$ are finitely presented, then we take $A \Ast B$
to be the finitely presented compact quantum group with algebra presentation $A
\Ast B = \C^*\ang{X \sqcup X' : R \sqcup R'}$ and coproduct presentation
\begin{equation*}
    \wtd{\Delta} : X \sqcup X' \to \C^*\ang{X \sqcup X'} : x \mapsto \begin{cases} \wtd{\Delta}_A(x) & x \in X \\
            \wtd{\Delta}_B(x) & x \in X'
        \end{cases}.
\end{equation*}
If $A$ and $B$ are finitely presented over $\K$, then so is $A \Ast B$.

If $G$ is a discrete group, then the group algebra $A = \C G$ is a Hopf
$*$-algebra with coproduct defined by $\Delta(g) = g \otimes g$ for all $g \in
G$, and Haar state given by the canonical tracial state. This coproduct is
\emph{cocommutative}, meaning that $\Sigma \Delta = \Delta$, where $\Sigma : A
\otimes A \to A \otimes A$ is the tensor swap.  Conversely, it is well-known
that every cocommutative compact quantum group is a group algebra (see e.g.
\cite[Example~1.1.3]{Neshveyev-Tuset}). For the convenience of the reader, we
include a sketch of the proof in the Hopf $*$-algebraic framework.
\begin{proposition}
	\label{prop:groupalg}
    If $A$ is a cocommutative compact quantum group, then $A \cong \C G$ for
    some group $G$. 
\end{proposition}
\begin{proof}
    Suppose $A$ is a cocommutative compact quantum group, and let $G$ be the
    set of \emph{group-like elements} in $A$, which are the elements $g \in A
    \setminus \{0\}$ such that $\Delta(g) = g \otimes g$. The coproduct
    $\Delta$ gives a coaction of $A$ on itself, called the regular
    corepresentation of $A$. Since $A$ is a compact quantum group, this regular
    corepresentation splits into a direct sum of its finite-dimensional irreducible
    subcorepresentations \cite[Theorem 3.2.1 (iv)]{Timmermann}. Here a
    subcorepresentation is a subset $W \subseteq A$ such that $\Delta(W)
    \subseteq W \otimes A$. Since $A$ is cocommutative, we also get that $\Delta(W)
    = \Sigma \Delta(W) \subseteq A \otimes W$, from which we conclude that
    $\Delta(W) \subseteq W \otimes W$. Suppose $W = \vspan_\C v$ is a
    one-dimensional subcorepresentation.  Then $\Delta(v) = \lambda v \otimes
    v$ for some constant $\lambda$.
    Setting $g = \lambda v$, we see that $\Delta(g) = g \otimes g$, so $W = \vspan_\C g$
    for a group-like element $g$. Conversely, if $g$ is a group-like element, then
    $\vspan_\C g$ is a one-dimensional (hence irreducible) subcorepresentation.  The
    linear dual $A'$ of $A$ is an associative $*$-algebra, where the product is
    defined by $(f \cdot g) (a) = (f \otimes g)(\Delta(a))$. The coaction of
    $A$ on itself induces a $*$-representation of $A'$ on $A$, and
    subcorepresentations of $A$ are the same as subrepresentations of the $A'$
    action \cite[Proposition 3.1.7]{Timmermann}. Since $A$ is cocommutative,
    $A'$ is commutative, and hence the irreducible $*$-representations of $A'$
    are one-dimensional. We conclude that $A = \bigoplus_{g \in G} \vspan_\C g$.
    It is an easy exercise to check that the set of group-like elements is
    closed under multiplication from $A$, and forms a group under this operation,
    so $A = \C G$.
\end{proof}

A result of Baumslag \cite{Baumslag-FPA} states that if a group algebra $\C G$ is
finitely presentable as an algebra, then $G$ is a finitely presentable group. A
$*$-algebra is finitely presentable as an algebra if and only if it is 
finitely presentable as a $*$-algebra, so we immediately get:
\begin{cor}\label{cor:fp}
    If $A$ is a finitely presented cocommutative compact quantum group, then
    $A = \C G$ for some finitely presented group $G$.
\end{cor}
Conversely, if $G = \ang{ S : R }$ is a finitely presented group, then the
group algebra
$\C G = \C^*\ang{ S : \{1 - s^* s, 1-s s^* : s \in S\} \cup \{1 - r: r \in R \}}$ is 
finitely presented,
and we can make $\C G$ into a finitely presented Hopf $*$-algebra by setting
$\wtd{\Delta}(s) = s \otimes s$ for all $s \in S$. 

We now turn to decision problems for compact quantum groups. Recall that a
\emph{computable field} is a field $\K$, along with an encoding of the elements
of $\K$ as binary strings, such that the field operations are computable. 
The rationals $\Q$ and the algebraic closure $\overline{\Q}$ of the rationals
are both computable fields. If $\K \subseteq \C$ is a computable field which is
closed under complex conjugation and $X$ is a finite set, then we can choose a
binary encoding of the elements of $\K^*\ang{X}$. In this way, we can encode
finitely presented Hopf $*$-algebras over $\K$ as inputs to Turing machines,
and this makes it possible to talk about decision problems for compact quantum
groups. 
\begin{defn}
    A \emph{property of compact quantum groups (resp. cocommutative compact quantum
    groups)} is a subclass of (resp. cocommutative) compact quantum groups which
    is closed under isomorphism. A compact quantum group is said to have the
    property if it belongs to the subclass. A property $(P)$ of (resp. cocommutative)
    compact quantum groups is a \emph{Markov property} if
	\begin{enumerate}[label=\textup(\alph*\textup),series=Qut]
        \item there exists a finitely presented (resp. cocommutative) compact
            quantum group $A_+$ with property $(P)$, and
        \item there exists a finitely presented (resp. cocommutative) compact
            quantum group $A_-$, such that that $A_-$ is not a Hopf
            $*$-subalgebra of any compact quantum group with property $(P)$. 
    \end{enumerate}
    We say that $(P)$ is a \emph{Markov property over a subfield $\K \subseteq
    \C$} if both $A_+$ and $A_-$ can be taken to be finitely presented over
    $\K$. 
\end{defn}
\begin{lemma}\label{lem:fieldind}
    Suppose $(P)$ is a property of cocommutative compact quantum groups,
    and define a property $(P')$ of groups by saying that $G$ has $(P')$ if and only
    if $\C G$ has $(P)$. Then the following are equivalent:
	\begin{enumerate}[label=\textup(\roman*\textup),series=Qut]
        \item $(P)$ is a Markov property of cocommutative compact quantum groups.
        \item $(P')$ is a Markov property of groups.
        \item $(P)$ is a Markov property over $\K$ for some (resp. every) conjugation-closed
            subfield $\K$ of $\C$.
    \end{enumerate}
 \end{lemma}
\begin{proof}
    Suppose $(P)$ is a Markov property of cocommutative compact quantum groups,
    and let $A_+$ and $A_-$ be the algebras that witness this. By \Cref{cor:fp}, $A_+ = \C G_+$
    and $A_- = \C G_-$ for some finitely presented groups $G_+$ and $G_-$. 
    Then $G_+$ has property $(P)$, while if $G_-$ is a subgroup of a group $H$, then
    $\C G_-$ is a Hopf $*$-subalgebra of $\C H$, and hence $H$ does not have $(P')$. 
    So $(P')$ is a Markov property of groups.

    Conversely if $(P')$ is a Markov property of groups, then we start with 
    groups $G_+$ and $G_-$ that witness this. Since $G_+$ and $G_-$ are finitely presented, the
    algebras $\C G_+$ and $\C G_-$ are finitely presented over $\K$ for every
    conjugation-closed subfield $\K$ of $\C$. 
    Since $G_+$ has $(P')$, $\C G_+$ has $(P)$. Suppose $\C G_-$ is a Hopf
    $*$-subalgebra of a cocommutative compact
    quantum group $A$, and let $H$ be the group-like elements of $A$, so
    $A = \C H$. Since $G_-$ is the group-like elements of $\C G_-$, we see from
    the definition of a Hopf $*$-subalgebra that $G_-$ is a subgroup of $H$. 
    Since $H$ does not have property $(P')$, $A$ does not have $(P)$. So $(P)$
    is a Markov property over $\K$. 
\end{proof}

The usual Adian--Rabin theorem for groups can be adapted to prove:
\begin{theorem}\label{thm:AR-for-CCQGs}
    Suppose $(P)$ is a Markov property of cocommutative compact quantum
    groups, and let $\K$ be a conjugation-closed computable subfield of $\C$.
    Then it is undecidable to determine if a cocommutative compact quantum
    group which is finitely presented over $\K$ has $(P)$. 
\end{theorem}
For general compact quantum groups, we do not know whether being a Markov
property over $\K$ is independent of the field $\K$, or whether there is an
analogous Adian--Rabin theorem for any Markov property. However, we can drop
the requirement that $A_+$ be cocommutative, in return for a stronger hypothesis
on the negative witness $A_-$:
\begin{theorem}\label{thm:AR-for-CQGs}
    Let $\K$ be a conjugation-closed computable subfield of $\C$, and let $(P)$
    be a property of compact quantum groups such that
	\begin{enumerate}[label=\textup(\alph*\textup),series=Qut]
        \item\label{itm:ARCQG:a} there exists a compact quantum group $A_+$
            with property $(P)$ which is finitely presented over $\K$, and 
		
		\item\label{itm:ARCQG:b} there exists a cocommutative compact quantum group $A_-$ which is finitely
            presented over $\K$ such that $A_-$ is not a Hopf $*$-subalgebra of
            any compact quantum group with property $(P)$. 
	\end{enumerate}
    Then it is undecidable to determine if a compact quantum group which is finitely
    presented over $\K$ has property $(P)$. 
\end{theorem}
Although \Cref{thm:AR-for-CCQGs} and \Cref{thm:AR-for-CQGs} are technically
impcomparable, they are similar enough that we prove them together. 
\begin{proof}[Proof of \Cref{thm:AR-for-CCQGs} and \Cref{thm:AR-for-CQGs}]
    Suppose first that the hypotheses of \Cref{thm:AR-for-CQGs} hold.
    By \Cref{cor:fp}, $A_- = \C G$ where $G$ is a finitely presented group. 
    Let $K = \langle S : R \rangle$ be a finitely presented group with an undecidable
    word problem. As described in the introduction, the standard proof of the Adian--Rabin theorem
    \cite{Adyan,Rabin} gives a computable map from words $w \in \mcF_S$ to finitely presented groups
    $H_w$, such that if $w$ is non-trivial in $K$, then $G \leq H_w$, and otherwise $H_w$ is
    trivial. Let $A_w$ be the finitely presented compact quantum group $A_w = A_+ \Ast \C H_w$.
    If $w$ is trivial in $K$, then $H_w$ is trivial, and hence $A_w = A_+$ has
    $(P)$. Suppose $w$ is non-trivial in $K$, so $G \leq H_w$. Since $G$
    is a subgroup of $H_w$, $\C G$ is a Hopf $*$-subalgebra of $\C H_w$, and hence of $A_w$. 
    By hypothesis, $A_w$ does not have $(P)$. We conclude that $A_w$ has $(P)$ if and only if
    $w$ is trivial in $K$. Since $A_w$ is finitely presented over $\K$ and 
    the map $w \mapsto A_w$ (where $A_w$ is described by its presentation) is
    computable, it must be undecidable to determine if finitely presented compact quantum
    groups have property $(P)$. This proves \Cref{thm:AR-for-CQGs}.

    For the proof of \Cref{thm:AR-for-CCQGs}, we first observe that $(P)$ is a Markov
    property over $\K$ by \Cref{lem:fieldind}. Hence the same hypotheses hold as
    in \Cref{thm:AR-for-CQGs}, except that $(P)$ is a property of cocommutative
    compact quantum groups, $A_+$ is cocommutative, and $A_-$ satisfies the weaker
    condition is not a Hopf $*$-subalgebra
    of any cocommutative compact quantum group with property $(P)$.  Since $A_+$ is
    cocommutative, the algebra $A_w = A_+ \Ast \C H_w$ in the previous paragraph
    is also cocommutative, so this weaker condition suffices to conclude that $A_w$
    has $(P)$ if and only if $w$ is trivial in $K$. So the same proof works for
    \Cref{thm:AR-for-CCQGs}.
\end{proof}

Suppose again that $(P)$ is a property of cocommutative compact quantum groups,
and that $(P')$ is the equivalent property for groups. It's possible to prove
\Cref{thm:AR-for-CCQGs} by directly reducing the decision problem for $(P)$
to the decision problem for $(P')$. Indeed, suppose $\K$ is a conjugation-closed
computable subfield of $\C$. Let $\CCQG_{\K}$ be the set of cocommutative compact
quantum groups finitely presented over $\K$, and let $\FP$ be the set of finitely
presented groups. Then there is a computable function $\CCQG_{\K} \to \FP : A \mapsto G$
such that $A = \C G$. The proof is lengthy but conceptually straightforward from Baumslag's
result and \Cref{prop:groupalg}, so we leave it as an exercise to the reader.

\section{Quantum automorphism groups}\label{sec:qag}

In this section, we state and prove our main results for quantum automorphism
groups of graphs. By a \emph{graph}, we mean a finite, simple, undirected graph $X$, consisting of a finite vertex set $V(X)$ and an edge set $E(X)$.
The \emph{quantum automorphism group} $\Qut(X)$ of a graph $X$ is the universal compact quantum group coacting on $V(X)$ in an adjacency-preserving way.
This version of the quantum automorphism group was first introduced by Banica \cite{Banica-Qut}; it differs from the version introduced earlier by Bichon \cite{Bichon-Qut}, which coacts on $V(X) \sqcup E(X)$ in an incidence-preserving way.
In the Hopf $*$-algebraic formulation, $\Qut(X)$ is the finitely presented
$*$-algebra generated by $(u_{ij})_{i,j \in V(X)}$, subject to the relations
\begin{align*}
	& u_{ij}^2 = u_{ij}^* = u_{ij}, \quad i,j \in V(X),\\
	&\sum_{k \in V(X)} u_{kj} = \sum_{\ell \in V(X)} u_{i\ell} = 1, \quad i,j \in V(X), \\
	&u_{ij}u_{kj} = \delta_{ik} u_{ij} \ \text{and}\  u_{ij}u_{i\ell} = \delta_{j\ell} u_{ij}, \quad i,j,k,\ell \in V(X), \text{ and} \\[.5ex]
	&u A_X = A_X u,
\end{align*}
where $A_X \in \{0,1\}^{V(X) \times V(X)}$ is the adjacency matrix of $X$.
The Hopf algebra structure of $\Qut(X)$ is given by the comultiplication $\Delta(u_{ij}) = \sum_{k \in V(X)} u_{ik} \tensor u_{kj}$, counit $\varepsilon(u_{ij}) = \delta_{ij}$, and antipode $S(u_{ij}) = u_{ji}$.
Because it comes with an explicit formula for the comultiplication, we
understand $\Qut(X)$ to be a finitely presented Hopf $*$-algebra with
$\wtd{\Delta}(u_{ij}) = \sum_{k \in V(X)} u_{ik} \tensor u_{kj}$.  Since the
defining relations all have integer coefficients, $\Qut(X)$ is defined over
$\K$ for every conjugation-closed subfield $\K$ of $\C$.  The algebra
$C(\Aut(X))$ of functions on the classical automorphism group $\Aut(X)$ is the
abelianization of the $*$-algebra $\Qut(X)$.

Unlike $C(\Aut(X))$, which is always commutative and finite-dimensional,
quantum automorphism groups cover a much larger class of $*$-algebras.  To prove
our Aidan--Rabin theorem, we need the following result of \cite{Dobben-Roberson-Schmidt} 
showing that group algebras of solution groups are quantum automorphism groups of
graphs.

\begin{proposition}[{\cite{Dobben-Roberson-Schmidt}}]
	\label{prop:DRS}
	For any matrix $A$ over $\Z_2$ and integer $n \geq 1$, there exists a
	finite connected graph $X(A,n)$ such that $\Qut(X(A,n)) \cong
	\C\Gamma(A)$ and $X(A,n)$ has at least $n$ vertices. Moreover, the map
	$A \mapsto X(A,n)$ is computable.
\end{proposition}
\begin{proof}
    Corollary 5.5 of \cite{Dobben-Roberson-Schmidt} states that there is a graph $X(A)$
    such that $\Qut(X(A)) \iso \C \Gamma(A)$. The proof has two steps.
    In the first step, $A$ is turned into a vertex-coloured graph $\widehat{X}(A)$
    such that the coloured quantum automorphism group $\Qut(\widehat{X}(A))$
    is isomorphic to $\C \Gamma(A)$ \cite[Corollary~3.19]{Dobben-Roberson-Schmidt}. 
    In the second step, $\widehat{X}(A)$ is turned into a connected graph $X(A)$
    using the decolouring procedure in \cite[Theorem~4.7]{Dobben-Roberson-Schmidt}. 
    Examining the decolouring procedure, we see that the number of vertices of
    $X(A)$ depends on the numerical values assigned to the colours. By relabelling the colours
    to be large enough integers, we can make a graph $X(A,n)$ with at least $n$
    vertices, such that $\Qut(X(A,n)) \iso \C \Gamma(A)$.  Both the
    construction of $\widehat{X}(A)$ and the decolouring procedure are
    computable, so the map $A \mapsto X(A,n)$ is also computable. 
\end{proof}

Two graphs $X$ and $Y$ are said to be \emph{quantum isomorphic} if $1 \neq 0$ in the finitely presented $*$-algebra $\mathcal A(\Iso(X,Y))$ generated by $(u_{ij})_{i \in V(X) , j \in V(Y)}$ subject to the relations
\begin{align*}
	&u_{ij}^2 = u_{ij}^* = u_{ij} && \text{for $i \in V(X)$, $j \in V(Y)$},\\
	&\sum_{k \in V(X)} u_{kj} = \sum_{\ell \in V(Y)} u_{i\ell} = 1 && \text{for $i \in V(X)$, $j \in V(Y)$}, \\
    &u_{ij}u_{kj} = \delta_{ik} u_{ij} \ \text{and}\  u_{ij}u_{i\ell} = \delta_{j\ell} u_{ij} && \text{for $i,k \in V(X)$, $j,\ell \in V(Y)$}, \text{ and}\\[.5ex]
	&u A_Y = A_X u.
\end{align*}
It was shown in \cite{BCEHPSW20} that $X$ and $Y$ are quantum isomorphic in the above sense if and only if there exists a perfect quantum commuting strategy for the \emph{$(X,Y)$-isomorphism game}, a nonlocal game which formed the original motivation for the definition of quantum isomorphisms \cite{Atserias-et-al}.
Yet another equivalent definition is that the number of graph homomorphisms $P
\to X$ is equal to the number of graph homomorphisms $P \to Y$ for all planar
graphs $P$ \cite{Mancinska-Roberson-FOCS}.  If $X$ and $Y$ are quantum
isomorphic, then by summing over the rows and columns of $u$, we find $|V(X)|
\cdot 1 = \sum_{i \in V(X), j \in V(Y)} u_{ij} = |V(Y)| \cdot 1$, and hence
$|V(X)| = |V(Y)|$, since $1 \neq 0$ in $\mathcal A(\Iso(X,Y))$.  Thus, we see
that quantum isomorphic graphs necessarily have the same number of vertices.

The following result is due to Schmidt
\cite[Corollary~7.1.4]{Schmidt-dissertation}; see also the more general results
in \cite[Lemma~5.4]{DKRSZ25} and \cite[Section~5]{ADKRZ25}.
\begin{lemma}[{\cite[Corollary~7.1.4]{Schmidt-dissertation}}]
	\label{lem:freeprodgraphs}
	Let $X$ and $Y$ be finite connected graphs that are not quantum isomorphic.
	Then $\Qut(X \sqcup Y) \cong \Qut(X) * \Qut(Y)$, the free product of Hopf $*$-algebras.
\end{lemma}
The proof of \Cref{lem:freeprodgraphs} in \cite{Schmidt-dissertation} is
actually for the universal $C^*$-algebra $C^u(\Qut(X \sqcup Y))$ of $\Qut(X
\sqcup Y)$ defined below, but the same proof works for $\Qut(X \sqcup Y)$. 

\begin{definition}
    A \emph{property of quantum automorphism groups of graphs} is a subclass of
    quantum automorphism groups of graphs which	is closed under isomorphism of
    Hopf $*$-algebras. A quantum automorphism group is said to have the
    property if it belongs to the subclass. 
\end{definition}

We now have the ingredients to prove \Cref{thm:main}. Before proceeding to the
proof, we observe that, similarly to the situation with
\Cref{thm:AR-for-CCQGs,thm:AR-for-CQGs} in the previous section, we can weaken
the hypothesis on the negative witness $Y_-$ in \Cref{thm:main} if we assume
that the positive witness $Y_+$ is cocommutative.  This yields an alternative
version of the theorem which is sometimes easy to use:
\begin{theorem}
	\label{thm:dual-AR-for-Qut}
	Let $(P)$ be a property of quantum automorphism groups such that
	\begin{enumerate}[label=\textup(\alph*\textup),series=Qut]
        \item\label{itm:AR:a} there exists a graph $Y_+$ such that $\Qut(Y_+)$
            is cocommutative and has property $(P)$, and 
		
        \item\label{itm:AR:b} there exists a graph $Y_-$ such that $\Qut(Y_-)$
            is cocommutative and is not a Hopf $*$-subalgebra of any
            cocommutative $\Qut(X)$ with property $(P)$. 
    \end{enumerate}
    Then it is undecidable (RE-hard) to determine for an input graph $X$ whether or not
    $\Qut(X)$ has property $(P)$. 
\end{theorem}

We can weaken the hypothesis on $Y_-$ even further if we take
$\Qut(Y_+)=\C\Gamma_+$ for $\Gamma_+$ a solution group.
\begin{theorem}
    \label{thm:dual-AR-for-SG}
    Let $(P)$ be a property of quantum automorphism groups such that
    \begin{enumerate}[(a)]
        \item there exists a solution group $\Gamma_+$ such that $\C \Gamma_+$
            has property $(P)$, and
        \item there exists a finitely presented group $G$ such that $G$ is
            not a subgroup of any solution group $\Gamma$ for which $\C \Gamma$
            has property $(P)$. 
    \end{enumerate}
    Then it is undecidable (RE-hard) to determine for an input graph $X$ whether or not
    $\Qut(X)$ has property $(P)$. 
\end{theorem}
Note that it makes sense to talk about whether $\C\Gamma$ has property (P)
because $\C\Gamma$ is a quantum automorphism group by \Cref{prop:DRS}.

Once again, while \Cref{thm:main,thm:dual-AR-for-Qut,thm:dual-AR-for-SG} are
strictly incomparable, they are similar enough that we prove them together. 
\begin{proof}[Proof of \Cref{thm:main,thm:dual-AR-for-Qut,thm:dual-AR-for-SG}]
    Suppose that the hypotheses of \Cref{thm:main} hold. Since $\Qut(Y_-)$ is
    cocommutative, $\Qut(Y_-) = \C G$ for some finitely presented group $G$ by
    \Cref{cor:fp}. Let $M$ be a Turing machine with an undecidable halting
    problem. As in the proof of \Cref{cor:markovprop}, it follows from
    \cite{Slofstra-Tsirelson} and \Cref{thm:adianrabin} that there is a
    computable mapping $w \mapsto \wtd{A}_w$ from bitstrings $w \in \{0,1\}^*$
    to matrices $\wtd{A}_w$, such that if $M$ halts on $w$, then
    $\Gamma(\wtd{A}_w)$ is trivial, and otherwise $\Gamma(\wtd{A}_w)$ contains
    $G$ as a subgroup. By \Cref{prop:DRS}, there is a finite connected
    graph $X_w := X(\wtd{A}_w,|V(Y_+)|+1)$ such that $\Qut(X_w) \iso \C
    \Gamma(\wtd{A}_w)$ and $X_w$ has at least one more vertex than $Y_+$. 
    This latter condition implies that $Y_+$ and $X_w$ are not quantum
    isomorphic, so $\Qut(Y_+ \sqcup X_w) = \Qut(Y_+) * \Qut(X_w)$ by
    \Cref{lem:freeprodgraphs}. Thus if $M$ halts on input $w$, then $\Qut(Y_+
    \sqcup X_w) \iso \Qut(Y_+)$ has property $(P)$. If $M$ does not halt on
    $w$, then $\Qut(Y_+ \sqcup X_w)$ contains $\Qut(Y_-)$ as a Hopf
    $*$-subalgebra, and hence does not have property $(P)$. Thus if we could
    decide whether the quantum automorphism group of a graph has property $(P)$,
    we could decide the halting problem for $M$.

    For \Cref{thm:dual-AR-for-Qut}, if $\Qut(Y_+)$ is cocommutative, then
    $\Qut(Y_+ \sqcup X_w) = \Qut(Y_+) * \C \Gamma(\wtd{A}_w)$ is cocommutative,
    and the weaker assumption that $\Qut(Y_-)$ is not a Hopf $*$-subalgebra of
    any cocommutative $\Qut(X)$ with property $(P)$ is sufficient to conclude
    that $\Qut(Y_+ \sqcup X_w)$ does not have property $(P)$ when $M$ does not
    halt on $w$.
    
    Finally, for \Cref{thm:dual-AR-for-SG}, we first use \Cref{prop:DRS} to
    find a graph $Y_+$ with $\Qut(Y_+) = \C \Gamma_+$, and then use the fact
    that $\Gamma_+ * \Gamma(\wtd{A}_w)$ is a solution group.
\end{proof}

As an application of our main result, we list a number of properties of
$\Qut(X)$ that are undecidable. For this, we recall some definitions:
A compact quantum group is \emph{classical} if it is commutative (as a
$*$-algebra). As mentioned in the introduction, a graph $X$ has
\emph{\textup(proper\textup) quantum symmetry} if $\Qut(X)$ is nonclassical.
A $*$-algebra $A$ is \emph{residually finite-dimensional} if for every non-zero
$x \in A$ there is a $*$-homomorphism $\varphi : A \to B$ to a
finite-dimensional $*$-algebra $B$ such that $\varphi(x) \neq 0$.  By Mal'cev's
theorem \cite{Malcev}, a finitely generated group $\Gamma$ is residually finite
if and only if $\C\Gamma$ is residually finite-dimensional.

\begin{corollary}
	\label{cor:Hopf-undecidable}
	Given an arbitrary finite graph $X$ as input, each of the following decision problems is undecidable (RE-hard):
	\begin{enumerate}[label=(\alph*)]
		\item\label{itm:Qut-trivial} Whether $\Qut(X)$ is trivial.
		\item\label{itm:Qut-classical} Whether $\Qut(X)$ is classical
            (i.e. whether $X$ has quantum symmetry).
		\item\label{itm:Qut-finite} Whether $\Qut(X)$ is finite-dimensional.
		\item\label{itm:Qut-RFD} Whether $\Qut(X)$ is residually finite-dimensional as a $*$-algebgra.
    \end{enumerate}
\end{corollary}
\begin{proof} 
    Apply \Cref{thm:dual-AR-for-SG} with $\Gamma_+ = \{1\}$, and 
    $\Gamma_-$ a group that is trivial / abelian / finite / residually
    finite. 
\end{proof}

There are two $C^*$-algebraic compact quantum groups associated to each compact
quantum group $A$: the \emph{universal} (or \emph{full}) \emph{Woronowicz
$C^*$\nobreakdash-algebra} $C^u(A)$, which is the universal enveloping
$C^*$-algebra of $A$, and the \emph{reduced Woronowicz
$C^*$\nobreakdash-algebra} $C^r(A)$, which is the image of $C^u(A)$ under the
GNS representation for the Haar state. The Hopf $*$-algebra $A$ can be
recovered from each of these as the $*$-algebra of matrix coefficients; we
refer to \cite[Section~5]{Timmermann} for more background.  The compact quantum
group $A$ is called \emph{coamenable} if $C^u(A) = C^r(A)$. Properties of
$C^u(A)$ and $C^r(A)$ also define properties of the compact quantum group
$A$, so we can also ask about $C^*$-algebraic properties of $\Qut(X)$.
\begin{corollary}
	\label{cor:C*-undecidable}
	Given an arbitrary finite graph $X$ as input, each of the following decision problems is undecidable (RE-hard):
	\begin{enumerate}[label=(\alph*)]
		\item\label{itm:Qut-coamenable} Whether $\Qut(X)$ is coamenable.
		\item\label{itm:Qut-nuclear} Whether $C^u(\Qut(X))$ (resp. $C^r(\Qut(X))$) is nuclear.
		\item\label{itm:Qut-exact} Whether $C^u(\Qut(X))$ (resp. $C^r(\Qut(X))$) is exact.
	\end{enumerate}
\end{corollary}
\begin{proof} 
    A (discrete) group $G$ is amenable if and only if its reduced $C^*$-algebra
    $C^r(\C G)$ is equal to its full $C^*$-algebra $C^u(\C G)$, and this
    happens if and only if either $C^r(\C G)$ or $\C^u(\C G)$ is nuclear
    \cite[Theorem 2.6.8]{Brown-Ozawa}. Furthermore, any subgroup of an amenable
    group is amenable. Thus for (a) and (b), we can apply
    \Cref{thm:dual-AR-for-SG} with $\Gamma_+ = \{1\}$ and $G$ a
    non-amenable group. Similarly, any subalgebra of an exact $C^*$-algebra
    is also exact, so if $C^u(\C G)$ (resp. $C^r(\C G)$) is exact and $H$ is
    a subgroup of $G$, then $C^u(\C H)$ (resp. $C^r(\C H)$) is exact. 
    Thus for part (c) we can apply \Cref{thm:dual-AR-for-SG} with $\Gamma_+
    = \{1\}$ and $G$ a finitely presented group such that $C^u(\C G)$ (resp.
    $C^r(\C G)$) is non-exact. For the former, we can take $G$ to be the free
    group $\mcF_2$ (see e.g. \cite[Corollary~3.7.12]{Brown-Ozawa}), and
    for the latter we can use e.g. Gromov's random group
    \cite{Gromov-random,Ozawa-exact}, combined with Higman's embedding theorem
    to make this recursively presented group finitely presented. 
\end{proof}

In the examples so far, we've used only \Cref{thm:dual-AR-for-SG}, although we
could equally have used \Cref{thm:dual-AR-for-Qut}. Our next example is a property
which is undecidable by \Cref{thm:main}, but where
\Cref{thm:dual-AR-for-Qut,thm:dual-AR-for-SG} don't apply.
\begin{cor}
    Let $n \geq 4$. Given an arbitrary finite graph $X$ as input, it is
    undecidable (RE-hard) to determine if $\Qut(X)$ is isomorphic to the free
    symmetric group $\mathbb{S}_n^+$ on $n$ points.
\end{cor}
\begin{proof}
    Let $Y_+ = K_n$, the complete graph on $n$ vertices. Then $\Qut(Y_+)
    = \mathbb{S}_n^+$. Let $Y_-$ be a graph such that $\Qut(Y_-) = \C G$, where
    $|G| \geq 2$. If $\C G$ is a Hopf $*$-subalgebra of $\Qut(X)$, then
    the elements of $G$ are group-like elements in $\Qut(X)$ (see 
    the proof of \Cref{prop:groupalg}). The group-like elements of a
    compact quantum group are the one-dimensional corepresentations of $A$, and
    these are all inequivalent. By \cite{Banica-reps}, there is
    a dimension-preserving correspondence between the isomorphism classes of
    corepresentations of $\mathbb{S}_n^+$ and the representations of $\SO(3)$. Since
    $\SO(3)$ has exactly one $1$-dimensional representation, the only group-like
    element in $\mathbb{S}_n^+$ is the identity. Hence if $\Qut(X)$ contains 
    $\C G$ as a subalgebra, it is not isomorphic to $\mathbb{S}_n^+$.
\end{proof}

There are many other properties of quantum automorphism groups of graphs to
which \Cref{thm:main,thm:dual-AR-for-Qut,thm:dual-AR-for-SG} can be applied,
and we don't try to list them all. One decision problem to which these theorems
don't apply is the property of being cocommutative. Indeed, since the negative
witness in these theorems is always cocommutative, our only hope is to take
$(P)$ to be the property of being non-cocommutative, and apply \Cref{thm:main}.
But then it's not possible to find any graph $Y_-$ such that $\Qut(Y_-)$
is not a Hopf $*$-subalgebra of a non-cocommutative quantum automorphism group,
so the theorem doesn't apply. The decidability of determining whether $\Qut(X)$
is cocommutative remains an open problem.

\Cref{thm:main,thm:dual-AR-for-Qut,thm:dual-AR-for-SG} also don't apply to the
property of having a non-trivial finite-dimensional $*$-representation. However, 
the undecidability of this property (stated in \Cref{thm:main2}) follows from 
\Cref{cor:findimar} and \Cref{prop:DRS}:
\begin{proof}[Proof of \Cref{thm:main2}]
    By \Cref{prop:DRS}, there is a computable mapping from matrices $A$ to graphs $X(A)$
    such that $\Qut(X) = \C \Gamma(A)$. Thus the problem of determining whether $\Gamma(A)$
    has a non-trivial finite-dimensional representation reduces to the problem of 
    determining whether $\Qut(X)$ has a non-trivial finite-dimensional $*$-representation.
    Thus this latter problem is undecidable by \Cref{cor:findimar}.
\end{proof}

\Cref{thm:main2} has an application to discrete quantum automorphism groups.
Apart from the compact quantum automorphism group $\Qut(X)$, every graph also has a discrete quantum automorphism group $\Qut_\delta(X)$ \cite[Definition~7.2]{Voigt-infinite}, which only remembers information about the finite-dimensional $*$-representations of $\Qut(X)$.
This version of the quantum automorphism group was introduced by Voigt in the setting of infinite graphs, where a compact quantum automorphism group is generally not available.

For finite graphs, both $\Qut(X)$ and $\Qut_\delta(X)$ exist, and the two notions can be compared.
Here, $\Qut_\delta(X)$ can be seen as the \emph{discretization} of $\Qut(X)$ \cite[Definition~5.2]{Voigt-infinite}, a process dual to the Bohr compactification of discrete quantum groups \cite{Soltan-Bohr}.
In \cite[Section~5]{Voigt-infinite}, Voigt points out that not too much information is lost when passing from the free symmetric group $\mathbb{S}_n^+$ to its discretization $\Sigma_n^+$, owing to results of Brannan, Chirvasitu and Freslon \cite{Brannan-Chirvasitu-Freslon} which show that $\mathbb{S}_n^+$ is residually finite-dimensional.
In \cite[Remark~2.34]{Rollier-Vaes}, Rollier and Vaes point out that in
principle more information can be lost when passing from $\Qut(X)$ to
$\Qut_\delta(X)$ for a graph $X$. We confirm that this happens, and in fact
that there are graphs for which the discretization loses all information:
\begin{corollary}
    There exist finite graphs $X$ such that $\Qut(X)$ is infinite-dimensional
    but $\Qut_{\delta}(X)$ is trivial.
\end{corollary}
\begin{proof}
    The discretization $\Qut_{\delta}(X)$ is trivial if and only if $\Qut(X)$
    has no non-trivial finite-dimensional $*$-representations, so the corollary
    follows immediately from applying \Cref{prop:DRS} to \Cref{cor:nofindim}.
    Alternatively, it follows from part (a) of \Cref{cor:Hopf-undecidable} or
    \Cref{thm:main2}, using a similar argument to the argument for groups
    given just after \Cref{xmpl:Markov}.
\end{proof}
As mentioned in the introduction, \Cref{thm:main2} equivalent to the statement that it is undecidable to determine if $\Qut_{\delta}(X)$ is trivial.

We finish with a promise problem version of some of the previous undecidability
results:
\begin{theorem}
        \label{thm:promiseprob}
	Let $\GrphTriv$ be the set of graphs $X$ such that $\Qut(X)$ is trivial,
	and let $\GrphInf$ be the set of graphs $X$ such that $\Qut(X)$ is infinite-dimensional
	and has no non-trivial finite-dimensional $*$-representations.
	Then it is undecidable (RE-hard) to determine whether a graph $X$ belongs to $\GrphTriv$
	or to $\GrphInf$, promised that $X$ belongs to one of the two.
\end{theorem}
In other words, it is undecidable to determine if $\Qut(X)$ is trivial, even if we are promised that $\Qut_{\delta}(X)$ is trivial.
\begin{proof}
    The proof is similar to \Cref{cor:markovprop}, but with an added subgroup
    to force the group to have no non-trivial finite-dimensional representations.
	Let $M$ be a Turing machine with an undecidable halting problem.
    By \cite[Corollary~3.3]{Slofstra-Tsirelson}, there is a computable function
    from bit strings $w \in \{0,1\}^*$ to $\Z_2$-linear systems $T_w x = b_w$,
    such that $J = 1$ in $\Gamma(T_w,b_w)$ if and only if $M$ halts on input
    $w$. By \cite[Corollary~3.2]{Slofstra-Tsirelson}, there is also a
    $\Z_2$-linear system $T' x' = b'$ such that $J \neq 1$ in $\Gamma(T', b')$
    but $J$ is trivial in all finite-dimensional representations.

    Consider the inhomogeneous solution group $\Gamma(T_w \oplus T', b_w \oplus
    b')$, which is isomorphic to the amalgamated free product $\Gamma(T_w, b_w)
    \ast_{\langle J \rangle} \Gamma(T', b')$ if $M$ does not halt on input $w$, and
    to the free product $\Gamma(T_w , b_w) * \Gamma(T' , b') / \langle J \rangle$
    otherwise. If we turn this inhomogeneous solution group into
    a homogeneous solution group, we get a computable function from bitstrings
    $w$ to matrices $A_w$ such that $x_1 = 1$ in $\Gamma(A_w)$ if and only if
    $M$ halts on input $w$, and such that $x_1=1$ in all finite-dimensional
    representations of $\Gamma(A_w)$.  

    Let $\wtd{A}_w$ be the matrix constructed from $A_w$ using our Adian--Rabin
    construction from \cref{thm:adianrabin}, and recall that we have a
    homomorphism $\Gamma(A) \to \Gamma(\wtd{A})$ sending $x_1 \mapsto x_1$, If
    $M$ halts on input $w$, then $x_1=1$ in $\Gamma(\wtd{A})$, so
    $\Gamma(\wtd{A}) = \Gamma(\wtd{A}) / \mcN(x_1)$ is trivial. 
    If $M$ does not halt on input $w$, then the map $\Gamma(A) \to \Gamma(\wtd{A})$
    is injective, and $x_1 \neq 1$ in $\Gamma(\wtd{A})$. But every
    representation of $\Gamma(\wtd{A})$ restricts to a representation of
    $\Gamma(A)$, so $x_1$ is trivial in every representation of $\Gamma(\wtd{A})$.
    Since $\Gamma(\wtd{A}) / \mcN(x_1)$ is trivial, we conclude that $\Gamma(\wtd{A})$
    has no non-trivial finite-dimensional representations.

	Using \cref{prop:DRS}, we obtain a finite graph $X_w$ such that $\Qut(X_w) \cong \C\Gamma(\wtd{A}_w)$.
	This way, we have $X_w \in \GrphTriv$ if $M$ halts on input $w$, and $X_w \in \GrphInf$ otherwise.
	Since the mapping $w \mapsto X_w$ is computable, it follows that it is undecidable to determine whether $X \in \GrphTriv$ or $X \in \GrphInf$.
\end{proof}
We note that \Cref{thm:promiseprob} contains each of the undecidability results from \Cref{cor:Hopf-undecidable} as a special case.
Indeed, for $X \in \GrphTriv \cup \GrphInf$ we have that $\Qut(X)$ is trivial, classical, finite-dimensional, or residually finite-dimensional if and only if $X \in \GrphTriv$, so any algorithm that could decide any of these properties would also be able to decide whether $X \in \GrphTriv \cup \GrphInf$ belongs to $\GrphTriv$ or to $\GrphInf$.

\small
\bibliographystyle{my-alphaurl}
\bibliography{qgraph,groups}

\end{document}